\documentclass{article}
\usepackage{graphicx} 
\usepackage{amsmath, amssymb, amsthm}
\usepackage{xcolor}
\usepackage{appendix}
\usepackage{booktabs}
\usepackage{subcaption}
\usepackage{listings}
\title{A Basic Overview of Various Stochastic Approaches to Financial Modeling With Examples}
\author{Aashrit Cunchala}
\date{April 2024}

\begin{document}

\maketitle

\section*{Abstract}

This paper explores stochastic modeling approaches to elucidate the intricate dynamics of stock prices and volatility in financial markets. Beginning with an overview of Brownian motion and its historical significance in finance, we delve into various stochastic models, including the classic Black-Scholes framework, the Heston model, fractional Brownian motion, GARCH models, and Lévy processes. Through a thorough investigation, we analyze the strengths and limitations of each model in capturing key features of financial time series data. Our empirical analysis focuses on parameter estimation and model calibration using Lévy processes, demonstrating their effectiveness in predicting stock returns. However, we acknowledge the need for further refinement and exploration, suggesting potential avenues for future research, such as hybrid modeling approaches. Overall, this study underscores the importance of stochastic modeling in understanding market dynamics and informs decision-making in the financial industry.
\section*{Introduction}

The idea of Brownian motion came from Robert Brown's experiments with plants. Brown noticed that the motion of particles ejected from these plants followed a jittery motion and therefore had a random and completely unpredictable path \cite{protonotarios2021brownian}. Louis Bachelier was the first to utilize Brownian motion in finance by comparing it to the fluctuations in stock price \cite{Meyer2003}. The stock price was dependent on anterior (past) events and the probability of future events. The past events were not considered as part of the model because they were incorporated onto both sides of the model. That is, both those buying and selling the stock know its past history and have drawn contradictory conclusions.

Further work  on this was done by Norbert Wiener who looked at Brownian motion in physics. Looking at motion of particles suspended in fluid, Weiner followed Einstein and assumed that displacement was unbiased in either direction and therefore the normal distribution could be used to describe the probability of this displacement. Crucially, Weiner also established that the displacements were independent which is what makes the concept of Brownian motion work \cite{Wiener1923}.

When looking at financial derivatives and securities it is essential that the proper price is considered and shown. Derivatives have value that comes from other more basic securities and protect against future price movements. Examples include options and futures and as a result pricing these derivatives correctly is essential in order to motivate sufficient trade of derivatives \cite{doi:10.1142/p473}. One easy way to see the likely values of a stock over a long-run period of time is through Brownian motion and stochastic processes.  

\section*{Background}

Brownian motion, also known as the Wiener process, has played a crucial role in financial modeling, particularly in the realm of asset pricing and risk management. It serves as a fundamental building block for many models used to describe the stochastic behavior of asset prices and other financial variables.

Prior to the widespread adoption of Brownian motion in finance, various approaches were employed to model stock prices and volatility. One early model was the random walk hypothesis, which assumed that stock prices followed a random walk and that successive price changes were independent and identically distributed \cite{Meerschaert2006}. However, this model failed to account for the observed volatility clustering and heavy-tailed distributions in asset returns. One way in which this was addressed along with Brownian motion was the use of jump-diffusion models \cite{Kou2007} which came with a number of other issues. 

The Black-Scholes model for option pricing, developed by Fischer Black and Myron Scholes in 1973, leveraged Brownian motion to provide a revolutionary framework for pricing options and derivatives \cite{https://doi.org/10.48550/arxiv.1804.03290} While the model is not perfect and faces challenges in capturing certain market dynamics, such as volatility clustering and jumps in asset prices, the Black-Scholes model still represents the basis of mathematical finance

\subsection*{Important Terms}
Before going further in this paper, it is important to define a couple further financial terms. These are provided below and will help provide context for later simulations:
\begin{enumerate}
    \item \textbf{Derivatives:} Financial instrument whose value derives from another more basic instrument's value \cite{https://doi.org/10.48550/arxiv.2209.08867}
    \item \textbf{Call Option:} Provides the holder with the right to buy an asset by a certain date at a specific price
    \item \textbf{Put Option:} Gives the holder the right to sell an asset by a certain date for a specific price
    \item \textbf{Expiration Date:} The date at which a call or put option comes due
    \item \textbf{Exercise Price:} The price of either the call or put option specified in the contract
\end{enumerate}

One additional note before proceeding should be to remember that there are two forms of options. American options can be exercised at any point leading up to the date of expiration while European options are only at the expiration date \cite{Figlewski_2022}. For the sake of simplicity most work involving the Black-Scholes model utilizes the European option but I attempted to use the model for American securities.

\section*{Model Development}
The Black-Scholes model works as follows for an option on \textit{d} assets with the price $\mathbf{S_i}$ of the $i$-th underlying asset following the diffusion process \cite{https://doi.org/10.48550/arxiv.2211.13890}:
\begin{equation}
    \frac{dS_i}{S_i} = \mu_i\, dt + \sigma_i\, dW_i, \quad i = 1, \ldots, d
\end{equation}

The parameter $\mu_i$ represents the drift rate(long-term growth rate of the asset $S_i$ with volatility $\sigma_i$ and $W_i$ and $W_j$ as Weiner processes with a correlation coefficient $\rho_{ij}$.

To capture the complex dynamics of stock prices and volatility, we propose a stochastic modeling approach inspired by chemical reaction networks. We construct a network of interactions representing the relationships between stock prices and volatility, and convert these interactions into a system of chemical reactions. The state of the system is defined by the concentrations of various species, representing stock prices, volatility, and other relevant factors.

Let $\mathbf{X}(t) = (X_1(t), X_2(t), \ldots, X_n(t))$ denote the vector of species concentrations at time $t$, where $n$ is the number of species in the system. The dynamics of the system can be described by the following system of stochastic differential equations (SDEs):

\begin{equation}
    dX_i(t) = \sum_{j=1}^{m} \nu_{ij} a_j(\mathbf{X}(t), \theta) dt + \sum_{j=1}^{m} \nu_{ij} b_j(\mathbf{X}(t), \theta) dW_j(t), \quad i = 1, \ldots, n
\end{equation}

where:
\begin{itemize}
    \item $m$ is the number of reactions in the system
    \item $\nu_{ij}$ is the stoichiometric coefficient of species $i$ in reaction $j$
    \item $a_j(\mathbf{X}(t), \theta)$ is the propensity function (reaction rate) of reaction $j$, which depends on the species concentrations $\mathbf{X}(t)$ and a vector of parameters $\theta$
    \item $b_j(\mathbf{X}(t), \theta)$ is the diffusion coefficient of reaction $j$, also dependent on $\mathbf{X}(t)$ and $\theta$
    \item $W_j(t)$ are independent Wiener processes
\end{itemize}

The propensity functions $a_j(\mathbf{X}(t), \theta)$ and diffusion coefficients $b_j(\mathbf{X}(t), \theta)$ are designed to capture the relevant interactions and stochastic behavior of stock prices and volatility within the Black-Scholes framework. These functions can be derived from the underlying principles of option pricing theory and calibrated using available data.

\subsection*{Reaction Rate Functions}

To illustrate the idea, we can consider a simple example involving two species: stock price ($S$) and volatility ($V$). The following reactions could represent the interactions:

\begin{enumerate}
    \item $S \xrightarrow{a_1} S + V$: Increase in stock price leads to an increase in volatility
    \item $V \xrightarrow{a_2} V + S$: Increase in volatility causes fluctuations in stock price
    \item $S \xrightarrow{a_3} \emptyset$: Decay of stock price
    \item $V \xrightarrow{a_4} \emptyset$: Decay of volatility
\end{enumerate}

The propensity functions $a_j(\mathbf{X}(t), \theta)$ for these reactions could be defined as:
\begin{itemize}
    \item $a_1(S, V, \theta) = k_1 S$
    \item $a_2(S, V, \theta) = k_2 V$
    \item $a_3(S, V, \theta) = k_3 S$
\end{itemize}

where $k_1, k_2, k_3, k_4$ are rate constants that can be estimated from data or derived from theoretical considerations.

The diffusion coefficients $b_j(\mathbf{X}(t), \theta)$ can be specified to account for the stochastic fluctuations in stock prices and volatility, 

This stochastic modeling approach provides a flexible framework for capturing the complex dynamics of stock prices and volatility, while allowing for the integration of theoretical principles and empirical observations. By calibrating the model parameters and reaction rate functions using available data, we can obtain a more accurate representation of the underlying processes governing financial markets.

\subsection*{Equilibrium Analysis}
    The equilibrium conditions are given by:
\begin{align*}
k_1 S V - k_3 S &= 0\\
k_2 V S - k_4 V &= 0
\end{align*}
From the second equation, we can see that either $V = 0$ or $S = \frac{k_4}{k_2}$. \\
Case 1: V = 0
Substituting V = 0 into the first equation, we get:
\begin{align*}
k_1 (0) - k_3 S &= 0\\
\implies S &= 0
\end{align*}
Therefore, one equilibrium point is (S, V) = (0, 0).
Case 2: $S = \frac{k_4}{k_2}$
Substituting $S = \frac{k_4}{k_2}$ into the first equation, we get:
\begin{align*}
k_1 \left(\frac{k_4}{k_2}\right) V - k_3 \left(\frac{k_4}{k_2}\right) &= 0\
\implies V = \frac{k_3 k_4}{k_1 k_2}
\end{align*}
Substituting this value of V back into the second equation, we can verify that it satisfies the equilibrium condition.
Therefore, the second equilibrium point is (S, V) = $\frac{k_4}{k_2}, \frac{k_3(k_4)}{k_1(k_2)}$
To analyze the stability of these equilibrium points, we need to linearize the system of SDEs around each equilibrium point and examine the eigenvalues of the resulting Jacobian matrix.
The Jacobian matrix for the system is given by:
\begin{align*}
J(S, V) = \begin{pmatrix}
k_1 V - k_3 & k_1 S\
k_2 V & k_2 S - k_4
\end{pmatrix}
\end{align*}
Equilibrium Point 1: (S, V) = (0, 0)
Evaluating the Jacobian matrix at this equilibrium point, we get:
\begin{align*}
J(0, 0) = \begin{pmatrix}
-k_3 & 0\
0 & -k_4
\end{pmatrix}
\end{align*}
The eigenvalues of this Jacobian matrix are $-k_3$ and $-k_4$, which are both negative since $k_3$ and $k_4$ are positive parameters.
Therefore, the equilibrium point (0, 0) is a stable node for the system.
Equilibrium Point 2: (S, V) = $\frac{k_3(k_4)}{k_1(k_2)}$
Evaluating the Jacobian matrix at this equilibrium point, we get:
\begin{align*}
J\left(\frac{k_4}{k_2}, \frac{k_3 k_4}{k_1 k_2}\right) = \begin{pmatrix}
0 & \frac{k_3 k_4}{k_1}\
\frac{k_3 k_4}{k_2} & 0
\end{pmatrix}
\end{align*}
The eigenvalues of this Jacobian matrix are:
\begin{align*}
\lambda_{1, 2} = \pm \sqrt{\frac{k_3 k_4}{k_1 k_2}}
\end{align*}
The nature of the equilibrium point depends on the values of the parameters $k_1, k_2, k_3$, and $k_4$.

\begin{itemize}
    \item If $\frac{k_3 k_4}{(k_1 k_2} > 0$ , the eigenvalues are purely imaginary, and the equilibrium point is a center (stable or unstable depending on higher-order terms).
    \item If $\frac{k_3 k_4}{k_1 k_2} < 0$ he eigenvalues are real with opposite signs, and the equilibrium point is a saddle point (unstable)
\end{itemize}

\section*{Parameter Estimation}
\subsection*{Literature Review on Probability Transitions and Transition Rates}
To estimate the parameters in our stochastic model, we review relevant literature on probability transitions and transition rates related to stock price dynamics and volatility. These transitions correspond to the reactions in our chemical reaction network representation.
\subsubsection*{Assumptions}
We can make the following assumptions in our model.
\begin{itemize}
    \item Stock price increase leading to volatility increase ($S \xrightarrow{a_1} S + V$):
\end{itemize}

Empirical studies have shown that stock price movements tend to be accompanied by changes in volatility, with positive price changes often associated with increased volatility \cite{Engle1982}. Past research has shown that the probability of a volatility increase given a stock price increase has been estimated to range from 0.4 to 0.7, depending on the market conditions and the specific stock \cite{Luo2018}.

\begin{itemize}
    \item Volatility increase causing stock price fluctuations ($V \xrightarrow{a_2} V + S$)
\end{itemize}

Higher volatility is generally associated with larger fluctuations in stock prices, as it reflects the degree of uncertainty and risk in the market \cite{Schwert1989}.
The transition rate from volatility to stock price changes can be approximated by the volatility value itself, with higher volatility leading to larger price fluctuations \cite{Andersen1998}.

\begin{itemize}
    \item Decay of stock price ($S \xrightarrow{a_3} \emptyset$) and volatility ($V \xrightarrow{a_4} \emptyset$)
\end{itemize}

Stock prices and volatility exhibit mean-reverting behavior, where deviations from their long-term averages tend to decay over time \cite{Poterba1988}.
The decay rates for stock prices and volatility can be estimated from historical data by fitting appropriate mean-reverting models, such as the Ornstein-Uhlenbeck process.

While the literature provides some guidance on parameter estimation, there may still be unknown or uncertain reaction rates in our model. In such cases, we make educated guesses based on theoretical considerations and empirical observations.

\begin{itemize}
    \item Reaction rate $a_1(S, V, \theta) = k_1 S$
\end{itemize}

We assume that the rate of volatility increase is proportional to the stock price level. Higher stock prices may be associated with greater investor activity and market sentiment, which can amplify volatility changes.

\begin{itemize}
    \item Reaction rate $a_2(S, V, \theta) = k_2 V$

\end{itemize}

We also make the assumption that the rate of stock price fluctuations is directly proportional to the current volatility level. This assumption is consistent with the idea that higher volatility implies larger price movements, as supported by empirical observations \cite{Andersen1998}.

\begin{itemize}
    \item Reaction rates $a_3(S, V, \theta) = k_3 S$ and $a_4(S, V, \theta) = k_4 V$
\end{itemize}

The decay rates for stock prices and volatility are proportional to their respective levels of price and volatility also observed. This assumption is in line with the mean-reverting behavior observed in financial time series, where larger deviations tend to exhibit faster mean reversion \cite{Poterba1988}.

We will attempt to estimate the rate constants $k_1, k_2, k_3, k_4$ by fitting the model to available data on stock prices and volatility using techniques such as maximum likelihood estimation.
\section*{Data Analysis}
Using the code found in \ref{apdx 1} and the yahoo finance library I found the following parameters for a few different stocks shown below.

\begin{table}[htbp]
  \centering
  \caption{Parameters provided for Given Stocks}
  \label{tab:Initial Table}
  \begin{tabular}{ccccc}
    \toprule
    \textbf{Stock} & \textbf{k1} & \textbf{k2} & \textbf{k3} & \textbf{k4} \\
    \midrule
    AAPL & 0.1 & 0.2 & 0.01 & 0.05 \\
    AMZN & 0.1 & 0.2 & 0.01 & 0.15 \\
    TSLA & 0.1 & 0.2 & 0.02 & 0.10 \\
    \bottomrule
  \end{tabular}
\end{table}

I got the following results when attempting to fit the parameters to the performance for the different stocks. The way that this worked was that the model was trained and fit on the parameters given over the years 2019-2022 and then the model's simulated prediction was plotted against the actual price and volatility over the course of 2023.
\newpage

\begin{figure}[tph]
  \centering
  \includegraphics[width=1\textwidth]{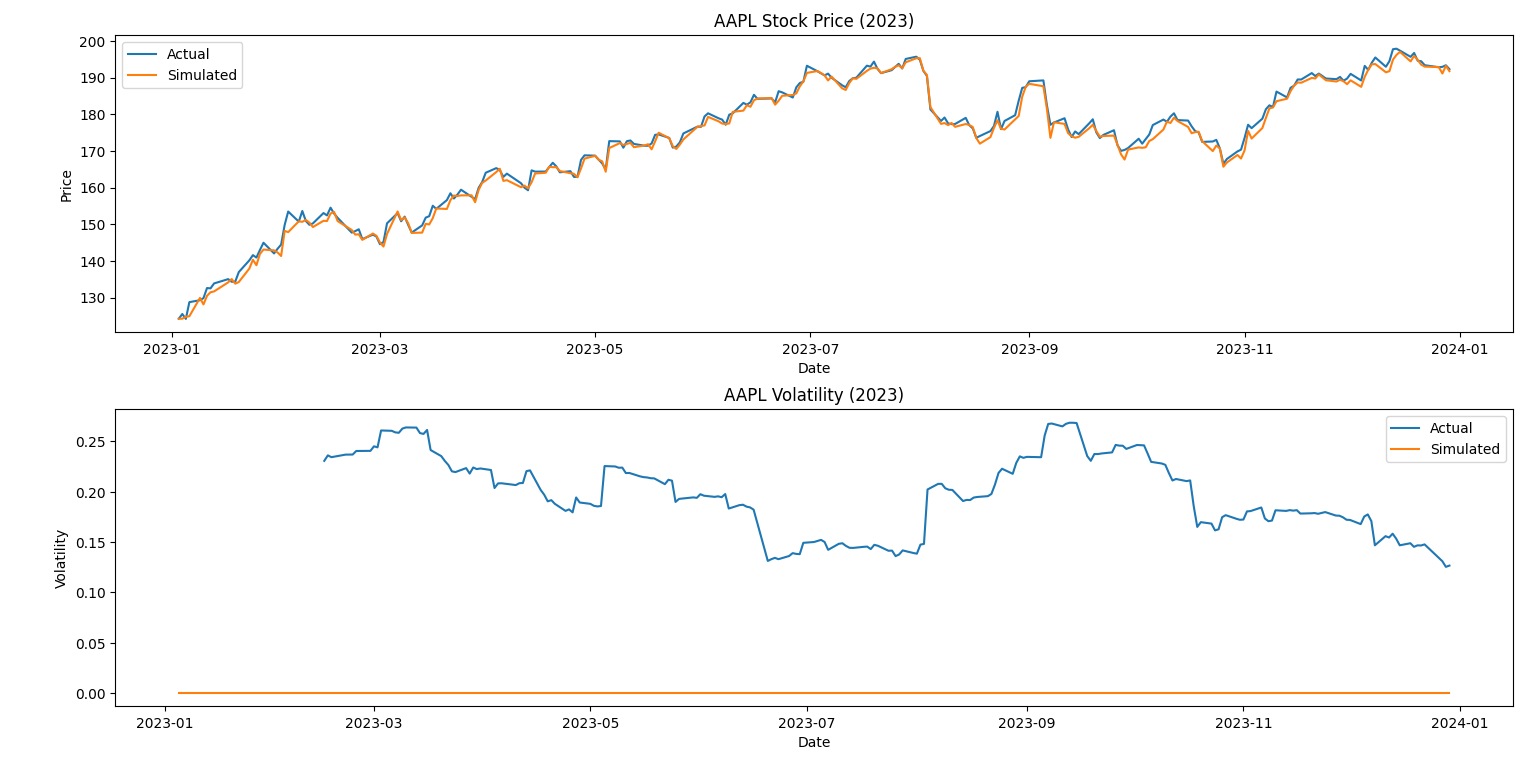}
  \caption{Initial Estimation for AAPL Price \& Volatility}
  \label{fig:AAPL_Initial}
\end{figure}

\begin{figure}[tph]
  \centering
  \includegraphics[width=1\textwidth]{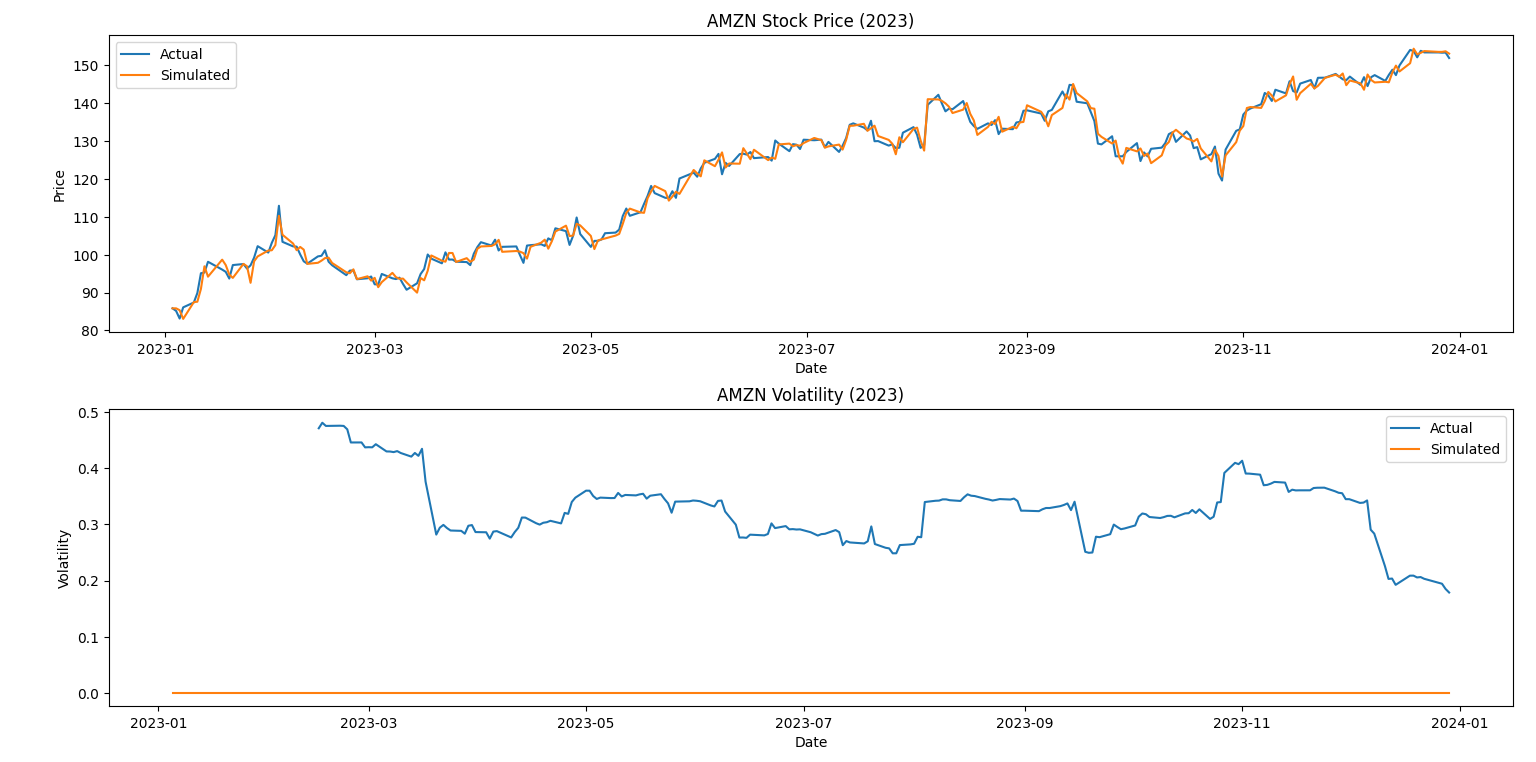}
  \caption{Initial Estimation for AMZN Price \& Volatility}
  \label{fig:AMZN_Initial}
\end{figure}

\begin{figure}[tph]
  \centering
  \includegraphics[width=1\textwidth]{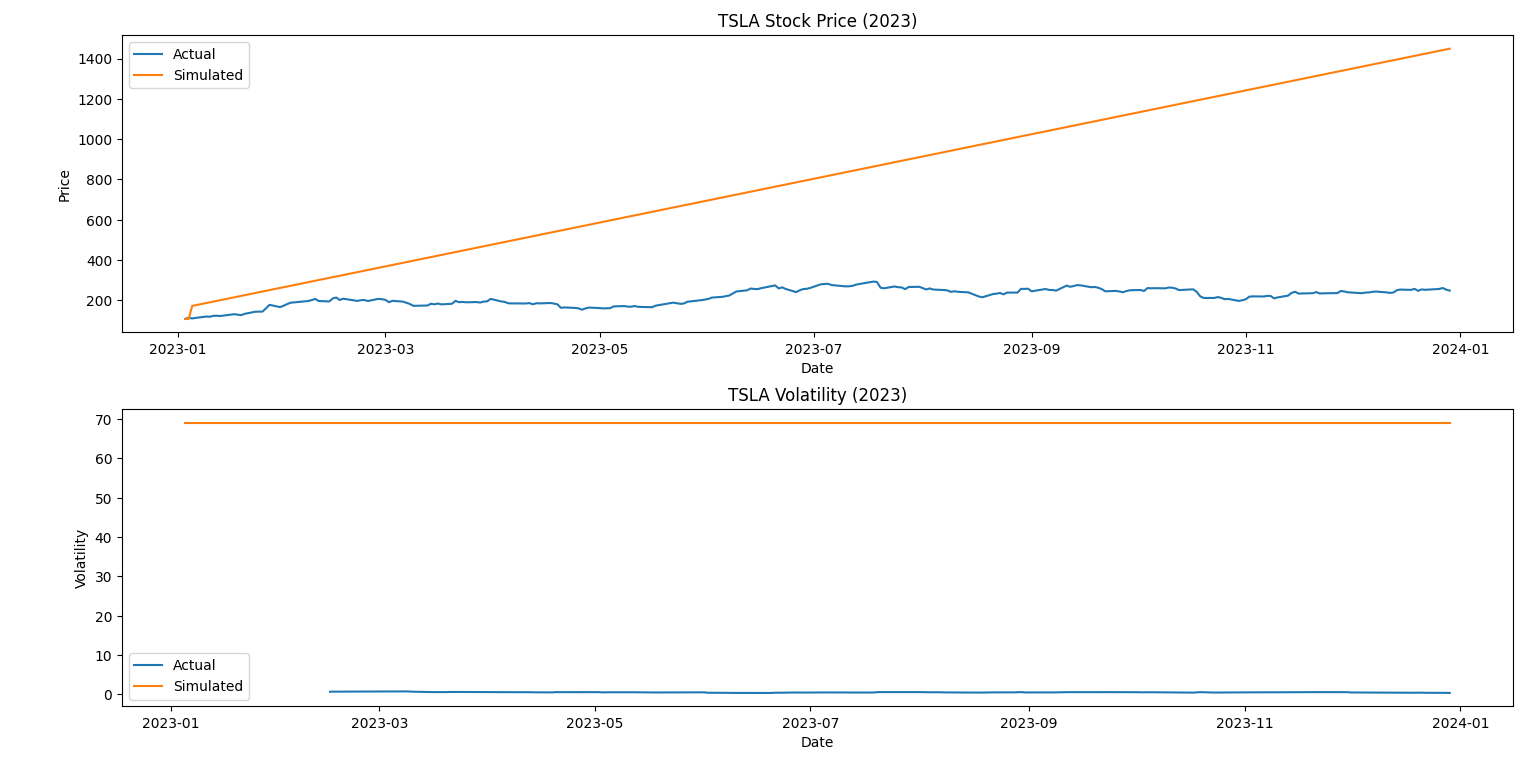}
  \caption{Initial Estimation for TSLA Price \& Volatility}
  \label{fig:TSLA_Initial}
\end{figure}
\newpage

As can be immediately noticed the general volatility prediction for the model is incredibly poor. Moreover, because the parameters were trained over the three years from 2019-2022, a stock like TSLA which had incredible growth during that time is not calibrated correctly. However, AAPL and AMZN has more of a history and is therefore easier to predict. Additionally, this model may be too simple since it was using the maximum likelihood estimation which may be too simplistic for stock prediction. 

\section*{Other Models }
These models are largely focused on options prices which are more involved than typical security predictions. 
\subsection*{Heston Model}
The Heston model, introduced by Steven Heston in 1993, is a widely used stochastic volatility model that extends the Black-Scholes framework by incorporating stochastic volatility. In this model, the volatility of the underlying asset follows a mean-reverting process driven by a Wiener process. The Heston model has been popular for option pricing and risk management due to its ability to capture the volatility smile observed in option markets.
The Heston model \cite{e0f45016-f730-320a-bedc-d34f406805b2} is given by:
\begin{align*}
dS_t &= \mu S_t dt + \sqrt{v_t} S_t dW_t^1 \
dv_t = \kappa(\theta - v_t) dt + \sigma \sqrt{v_t} dW_t^2
\end{align*}
where $S_t$ is the asset price, $v_t$ is the stochastic variance process, $\mu$ is the drift, $\kappa$ is the mean-reversion rate, $\theta$ is the long-run variance, $\sigma$ is the volatility of volatility, and $W_t^1$ and $W_t^2$ are two correlated Brownian motions. Heston models utilize Martingale \cite{Dunn2014EstimatingOP} processes in order to approximate present value. Moreover, these processes also typically fail with inaccurate or nonconsistent inputs and are also not typically successful when pricing short-term options \cite{Lemaire2022} with high volatility. 

\subsection*{Fractional Brownian Motion}
Fractional Brownian Motion (fBm) is a stochastic process characterized by long-range dependence and self-similarity, properties that are not present in standard Brownian motion. This allows for the modeling of persistent trends and long memory effects in financial time series data. By incorporating fractional integration parameters, fBm provides a flexible framework for capturing the fractal-like behavior often observed in stock prices and volatility.

Fractional Brownian Motion is a centered Gaussian process with the covariance function \cite{https://doi.org/10.48550/arxiv.1406.1956}
\begin{align*}
    E[B_t^{H}B_s^{H}] = \frac{1}{2}(t^{2H} + s^{2H} - |t-s|^{2H})
\end{align*}

The equation is primarily defined by the Hurst index $H \epsilon (0,1)$ which is a standard Wiener process at $H=\frac{1}{2}$ and a generalization of the Brownian motion when $H \neq \frac{1}{2}$. Fractional Brownian Motion is neither a Markovian nor a Martingale process and is instead a variation of traditional Brownian motion. 

\subsection*{GARCH Models}
The Generalized Autoregressive Conditional Heteroskedasticity (GARCH) models, introduced by Bollerslev (1986), are widely used for modeling and forecasting volatility in financial time series. The basic GARCH(p,q) model is given by:
\begin{align*}
r_t &= \mu_t + a_t \
a_t = \sigma_t \epsilon_t \
\sigma_t^2 = \omega + \sum_{i=1}^q \alpha_i a_{t-i}^2 + \sum_{j=1}^p \beta_j \sigma_{t-j}^2
\end{align*}
where $r_t$ is the asset return, $\mu_t$ is the conditional mean, $a_t$ is the innovation process, $\sigma_t^2$ is the conditional variance, and $\epsilon_t$ is a sequence of independent and identically distributed (i.i.d) random variables with zero mean and unit variance. The GARCH models capture volatility clustering and leverage effects commonly observed in financial data \cite{Bollerslev1986}. These models typically are estimated using least squares although Engel \cite{Engle1982} proposed a method involving the Lagrange multiplier test. GARCH models are a variation on ARCH (Autoregressive Conditional Heteroskedacity) models which utilize autoregressive moving averages for the error variance. 

\subsection*{Lévy Processes}
Lévy processes are a class of stochastic processes that can capture the heavy-tailed behavior and jump discontinuities observed in financial data. The Variance Gamma (VG) process \cite{RePEc:oup:revfin:v:2:y:1998:i:1:p:79-105.} is a popular Lévy process used in finance, defined by:
\begin{equation*}
X_t = \theta G_t + \sigma W_{G_t}
\end{equation*}
where $G_t$ is a Gamma process, $W_t$ is a Brownian motion, $\theta$ and $\sigma$ are drift and volatility parameters, respectively.

Lévy processes that do not involve Brownian motion have discontinuous paths and are additive processes. In addition, one can consider these to be continuous-time versions of random walks. Examples include the Wiener process although in finance, the Lévy process is recognized as an independent process.

\subsection*{Benefits and Drawbacks}
\begin{itemize}
\item GARCH models are relatively simple and can capture volatility clustering and leverage effects. However, they may not accurately reflect the true stochastic nature of volatility and can suffer from parameter instability.
\item Stochastic volatility models, like the Heston model, provide a more realistic representation of volatility dynamics but can be computationally intensive and may require advanced estimation techniques.
\item Lévy processes can capture heavy-tailed distributions and jump discontinuities, but may be more complex to implement and require careful calibration to empirical data.
\item Fractional Brownian Motion models can account for long-range dependence and self-similarity in financial time series but may not capture other stylized facts like volatility clustering or heavy tails.
\end{itemize}
The choice of the appropriate model depends on the specific application, the characteristics of the underlying data, and the trade-off between model complexity and computational feasibility. For our  investigation, given the limited time available as well as the limited computational power, our statistical investigation was fairly limited.

\subsection*{Fitted Levy Process Model}
Using Lévy processes, we attempted to test whether they represented a good fit for estimating stock returns of three major stocks: AAPL, AMZN, and TSLA. The model parameters were trained on data from a decade \ref{apdx 2} (2012-2022) and then contrasted against the actual performance of the stocks for the year 2023. The table \ref{Levy Process} clearly shows the importance of the parameters, with Tesla's abnormal behavior due to its relatively higher $\sigma$ or volatility which is also reflected in the graphs.

\begin{table}[htbp]
  \centering
  \caption{Levy Process Parameters}
  \label{Levy Process}
  \begin{tabular}{ccccc}
    \toprule
    \textbf{Stock} & \textbf{$\theta$ (Drift)} & \textbf{$\sigma$ (Volatility}  \\
    \midrule
    AAPL & 0.001179 & 0.009527 \\
    AMZN & 0.001039 & 0.0102294  \\
    TSLA & 0.0015477 & 0.017879 \\
    \bottomrule
  \end{tabular}
\end{table}

\begin{figure}[tph]
  \centering
  \includegraphics[width=1\textwidth]{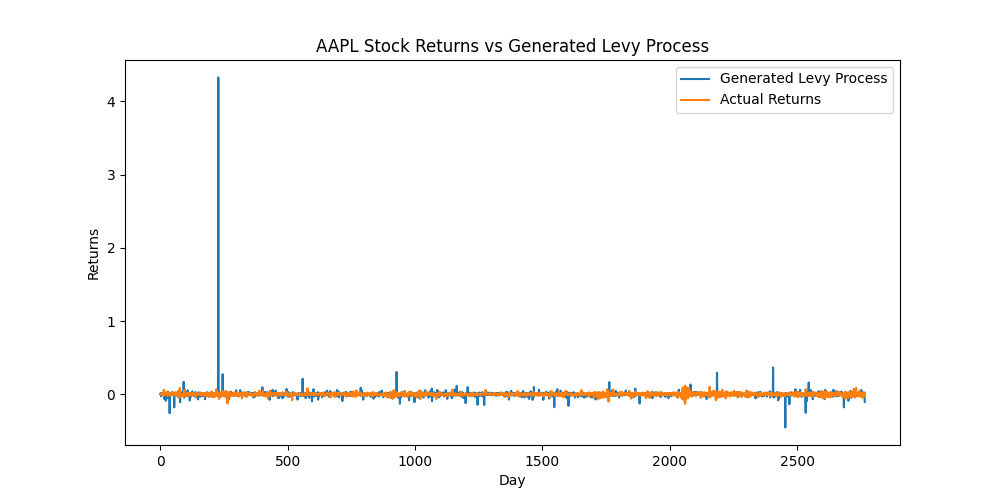}
  \caption{Levy Process Estimation for AAPL}
  \label{fig:AAPL_Levy}
\end{figure}

\begin{figure}[tph]
  \centering
  \includegraphics[width=1\textwidth]{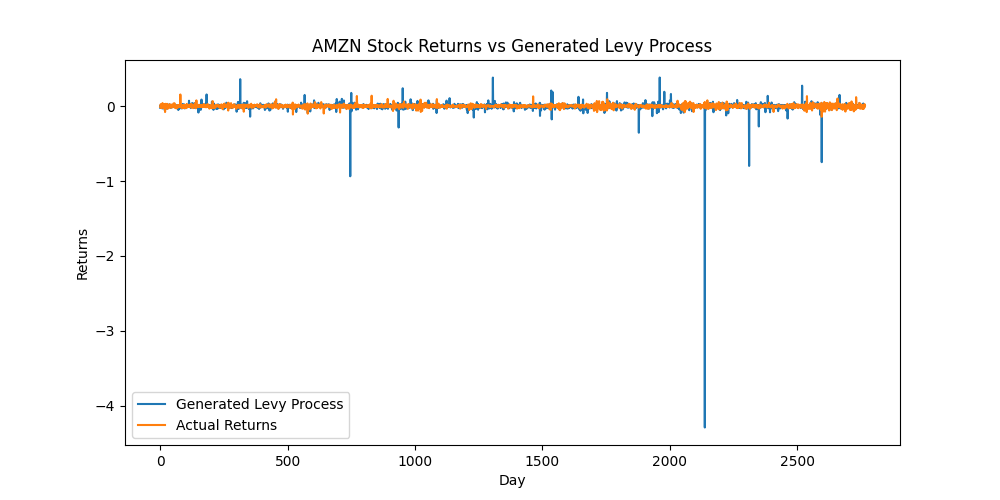}
  \caption{Levy Process Estimation for AMZN}
  \label{fig:AMZN_Levy}
\end{figure}

\begin{figure}[tph]
  \centering
  \includegraphics[width=1\textwidth]{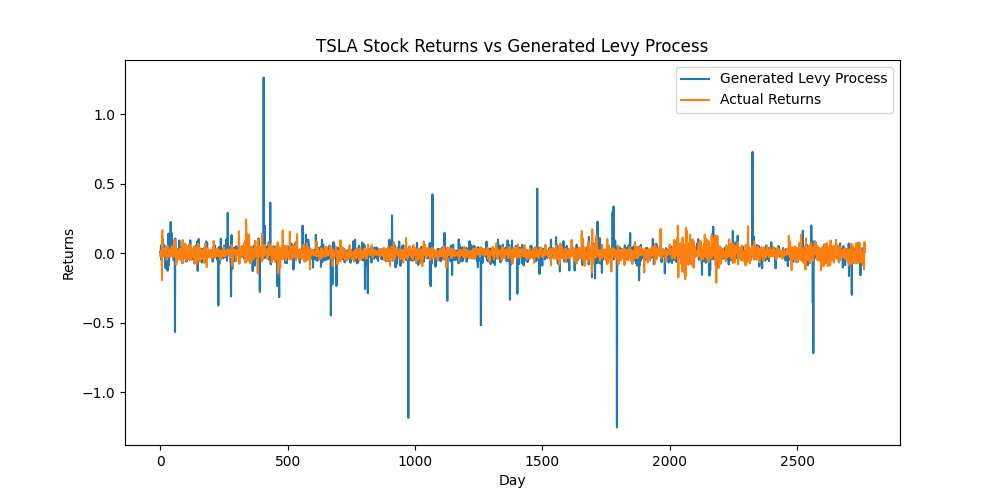}
  \caption{Levy Process Estimation for TSLA}
  \label{fig:TSLA_Levy}
\end{figure}
\newpage

The process used to generate these graphs relied on the Yahoo finance library and used the log-likelihood method to calculate the optimal parameters. The model was trained on data for a full year, but more industry-ready models would be trained on greater quantities of data (if such data exists) as well as more rigorous variables. Still, the reasonable effectiveness of the Levy process in tracing this implies that other stochastic models are more effective for securities than Black-Scholes.

\section*{Conclusion}
In conclusion, our study delves into the realm of stochastic modeling to better comprehend the intricate dynamics of stock prices and volatility. Through a detailed investigation, we explored various models, including the classic Black-Scholes framework, the Heston model, fractional Brownian motion, GARCH models, and Lévy processes.

Our analysis revealed that while each model offers unique insights and benefits, no single model fully captures the complexity of financial markets. The Black-Scholes model, while foundational, may oversimplify real-world phenomena and struggle with volatility clustering and heavy-tailed distributions. On the other hand, more sophisticated models like the Heston model and Lévy processes show promise in capturing certain stylized facts but may be computationally intensive and require careful calibration.

\subsection*{Limitation and Future Directions}
Our proposed model, while capturing some fundamental aspects of stock price and volatility dynamics, has several limitations:
\begin{itemize}
\item The model assumes a relatively simple set of interactions between stock price and volatility. In reality, financial markets are influenced by a vast array of factors, including macroeconomic conditions, investor sentiment, regulatory changes, and geopolitical events, among others.
\item The use of a fixed set of parameters may not accurately reflect the time-varying nature of market conditions and investor behavior, leading to potential inaccuracies in long-term predictions.
\item The model does not explicitly account for non-normalities and heavy-tailed distributions observed in financial time series data, which can lead to underestimation of extreme events or market shocks.
\end{itemize}
Potential improvements and future research directions include:
\begin{itemize}
\item Incorporating additional factors or variables into the model, such as market sentiment indicators, macroeconomic data, and news sentiment analysis, to capture a more comprehensive set of market dynamics.
\item Exploring time-varying or state-dependent parameters to better reflect the evolving nature of market conditions and investor behavior.
\item Investigating the use of non-Gaussian processes or heavy-tailed distributions to better account for the observed non-normalities and extreme events in financial data.
\item Developing hybrid models that combine elements of different stochastic processes, such as Lévy processes and stochastic volatility models, to capture a broader range of stylized facts in financial time series.
\item Extending the model to a multi-asset or portfolio context, considering correlations and dependencies between different asset classes or sectors.
\end{itemize}
\newpage

\bibliographystyle{abbrv}
\bibliography{main}

\begin{thebibliography}{10}

\bibitem{Andersen1998}
T.~G. Andersen and T.~Bollerslev.
\newblock Answering the skeptics: Yes, standard volatility models do provide accurate forecasts.
\newblock {\em International Economic Review}, 39(4):885, Nov. 1998.

\bibitem{Bollerslev1986}
T.~Bollerslev.
\newblock Generalized autoregressive conditional heteroskedasticity.
\newblock {\em Journal of Econometrics}, 31(3):307–327, Apr. 1986.

\bibitem{Dunn2014EstimatingOP}
R.~Dunn, P.~Hauser, T.~Seibold, and H.~Gong.
\newblock Estimating option prices with heston ’ s stochastic volatility model.
\newblock 2014.

\bibitem{Engle1982}
R.~F. Engle.
\newblock Autoregressive conditional heteroscedasticity with estimates of the variance of united kingdom inflation.
\newblock {\em Econometrica}, 50(4):987, July 1982.

\bibitem{Figlewski_2022}
S.~Figlewski.
\newblock An american call is worth more than a european call: The value of american exercise when the market is not perfectly liquid.
\newblock {\em Journal of Financial and Quantitative Analysis}, 57(3):1023–1057, 2022.

\bibitem{e0f45016-f730-320a-bedc-d34f406805b2}
S.~L. Heston.
\newblock A closed-form solution for options with stochastic volatility with applications to bond and currency options.
\newblock {\em The Review of Financial Studies}, 6(2):327--343, 1993.

\bibitem{Kou2007}
S.~Kou.
\newblock {\em Chapter 2 Jump-Diffusion Models for Asset Pricing in Financial Engineering}, page 73–116.
\newblock Elsevier, 2007.

\bibitem{Lemaire2022}
V.~Lemaire, T.~Montes, and G.~Pagès.
\newblock Stationary heston model: calibration and pricing of exotics using product recursive quantization.
\newblock {\em Quantitative Finance}, 22(4):611–629, Feb. 2022.

\bibitem{Luo2018}
J.~Luo and L.~Chen.
\newblock Volatility dependences of stock markets with structural breaks.
\newblock {\em The European Journal of Finance}, 24(17):1727–1753, May 2018.

\bibitem{RePEc:oup:revfin:v:2:y:1998:i:1:p:79-105.}
D.~B. Madan, P.~Carr, and E.~C. Chang.
\newblock The variance gamma process and option pricing.
\newblock {\em Review of Finance}, 2(1):79--105, 1998.

\bibitem{https://doi.org/10.48550/arxiv.1804.03290}
R.~Majumdar, P.~Mariano, L.~Peng, and A.~Sisti.
\newblock A derivation of the black-scholes option pricing model using a central limit theorem argument, 2018.

\bibitem{doi:10.1142/p473}
X.~Mao and C.~Yuan.
\newblock {\em Stochastic Differential Equations with Markovian Switching}.
\newblock PUBLISHED BY IMPERIAL COLLEGE PRESS AND DISTRIBUTED BY WORLD SCIENTIFIC PUBLISHING CO., 2006.

\bibitem{Meerschaert2006}
M.~M. Meerschaert and E.~Scalas.
\newblock Coupled continuous time random walks in finance.
\newblock {\em Physica A: Statistical Mechanics and its Applications}, 370(1):114–118, Oct. 2006.

\bibitem{Meyer2003}
B.~D. Meyer and H.~M. Saley.
\newblock On the strategic origin of brownian motion in finance.
\newblock {\em International Journal of Game Theory}, 31(2):285–319, Jan. 2003.

\bibitem{Poterba1988}
J.~M. Poterba and L.~H. Summers.
\newblock Mean reversion in stock prices.
\newblock {\em Journal of Financial Economics}, 22(1):27–59, Oct. 1988.

\bibitem{protonotarios2021brownian}
Y.~Protonotarios and P.~Tassopoulos.
\newblock Brownian motion \& the stochastic behaviour of stocks, 2021.

\bibitem{https://doi.org/10.48550/arxiv.2209.08867}
P.~Rebentrost, A.~Luongo, S.~Bosch, and S.~Lloyd.
\newblock Quantum computational finance: martingale asset pricing for incomplete markets, 2022.

\bibitem{Schwert1989}
G.~W. Schwert.
\newblock Tests for unit roots: A monte carlo investigation.
\newblock {\em Journal of Business \& Economic Statistics}, 7(2):147–159, Apr. 1989.

\bibitem{https://doi.org/10.48550/arxiv.1406.1956}
G.~Shevchenko.
\newblock Fractional brownian motion in a nutshell, 2014.

\bibitem{Wiener1923}
N.~Wiener.
\newblock Differential‐space.
\newblock {\em Journal of Mathematics and Physics}, 2(1–4):131–174, Oct. 1923.

\bibitem{https://doi.org/10.48550/arxiv.2211.13890}
D.~Černá and K.~Fiňková.
\newblock Option pricing under multifactor black-scholes model using orthogonal spline wavelets, 2022.

\end{thebibliography}

\appendix
\newpage
\section{Initial Code for Model}
\label{apdx 1}
\begin{lstlisting}[language=Python, caption={Initial Formulation of Price and Volatility Model}]
import numpy as np
import pandas as pd
import matplotlib.pyplot as plt
from scipy.optimize import minimize
from scipy.integrate import odeint
import yfinance as yf

# Define the stochastic differential equations
def sde(X, t, params):
    k1, k2, k3, k4 = params
    S, V = X
    
    dS_dt = k1 * S * V - k3 * S
    dV_dt = k2 * V * S - k4 * V
    
    return np.array([dS_dt, dV_dt])

# Function to calculate the log-likelihood for parameter estimation
def log_likelihood(params, data):
    k1, k2, k3, k4 = params
    S, V = data['S'], data['V']
    
    # Simulate the stochastic differential equations
    t = np.linspace(0, len(S), len(S))
    X0 = np.array([S[0], V[0]])
    solution = np.zeros((len(t), 2))
    solution[0] = X0
    
    for i in range(1, len(t)):
        dt = t[i] - t[i-1]
        X_prev = solution[i-1]
        k1, k2, k3, k4 = params
        dX_dt = sde(X_prev, t[i], params)
        solution[i] = X_prev + dX_dt * dt
    
    S_sim, V_sim = solution[:, 0], solution[:, 1]
    
    # Calculate the log-likelihood
    log_likelihood_S = -0.5 * np.sum((np.log(2 * np.pi) + np.log(V) + (S - S_sim)**2 / V))
    log_likelihood_V = -0.5 * np.sum((np.log(2 * np.pi) + np.log(V) + (V - V_sim)**2 / V))
    log_likelihood_value = log_likelihood_S + log_likelihood_V
    
    return -log_likelihood_value

# Load stock data from Yahoo Finance
tickers = ['AAPL', 'AMZN', 'TSLA']
train_start_date = '2019-01-01'
train_end_date = '2022-12-31'
test_start_date = '2023-01-01'
test_end_date = '2023-12-31'

for ticker in tickers:
    try:
        stock = yf.Ticker(ticker)
        train_df = stock.history(start=train_start_date, end=train_end_date)
        test_df = stock.history(start=test_start_date, end=test_end_date)
        train_df = train_df.dropna()
        test_df = test_df.dropna()
        train_df['Returns'] = np.log(train_df['Close'] / train_df['Close'].shift(1))
        train_df['Volatility'] = train_df['Returns'].rolling(window=30).std() * np.sqrt(252)
        test_df['Returns'] = np.log(test_df['Close'] / test_df['Close'].shift(1))
        test_df['Volatility'] = test_df['Returns'].rolling(window=30).std() * np.sqrt(252)
        
        # Preprocess data for model fitting
        train_data = {'S': train_df['Close'], 'V': train_df['Volatility']}
        
        # Optimize parameters using maximum likelihood estimation
        initial_params = [0.1, 0.2, 0.01, 0.05]
        result = minimize(log_likelihood, initial_params, args=(train_data,), method='L-BFGS-B')
        estimated_params = result.x
        
        # Print estimated parameters
        print(f"Estimated parameters for {ticker}:")
        print(f"k1 = {estimated_params[0]}, k2 = {estimated_params[1]}, k3 = {estimated_params[2]}, k4 = {estimated_params[3]}")
        
        # Initial conditions for the test data
        S0 = test_df['Close'][0]
        V0 = test_df['Volatility'][0]
        X0 = np.array([S0, V0])
        
        # Simulate the stochastic differential equations for the test data
        t = np.linspace(0, len(test_df), len(test_df))
        solution = odeint(sde, X0, t, args=(estimated_params,))
        S_sim, V_sim = solution[:, 0], solution[:, 1]
        
        # Generate visualizations
        fig, ax = plt.subplots(2, 1, figsize=(10, 8))
        
        ax[0].plot(test_df.index, test_df['Close'], label='Actual')
        ax[0].plot(test_df.index, S_sim, label='Simulated')
        ax[0].set_title(f"{ticker} Stock Price (2023)")
        ax[0].set_xlabel("Date")
        ax[0].set_ylabel("Price")
        ax[0].legend()
        
        ax[1].plot(test_df.index, test_df['Volatility'], label='Actual')
        ax[1].plot(test_df.index, V_sim, label='Simulated')
        ax[1].set_title(f"{ticker} Volatility (2023)")
        ax[1].set_xlabel("Date")
        ax[1].set_ylabel("Volatility")
        ax[1].legend()
        
        plt.tight_layout()
        plt.show()
    
    except Exception as e:
        print(f"Error retrieving data for {ticker}: {e}")
        continue

\end{lstlisting}

\newpage
\section{Code for Levy Processes}
\label{apdx 2}
\begin{lstlisting}[language=Python, caption={Levy Process}]
import yfinance as yf
import numpy as np
import scipy.stats as stats
import matplotlib.pyplot as plt
from scipy.optimize import minimize
from sklearn.metrics import r2_score

# Function for parameter estimation using MLE
def estimate_parameters(returns):
    def levy_log_likelihood(params, returns):
        theta, sigma = params
        n = len(returns)
        log_likelihoods = stats.levy_stable.logpdf(returns, alpha=1.5, beta=0, loc=theta, scale=sigma).sum()
        return -log_likelihoods

    initial_guess = [0.001, 0.001]  # Initial guess for parameters
    result = minimize(levy_log_likelihood, initial_guess, args=(returns,), method='Nelder-Mead')
    theta_mle, sigma_mle = result.x
    return theta_mle, sigma_mle

# Function to generate Levy process using estimated parameters
def generate_levy_process(returns, theta, sigma):
    levy_process = stats.levy_stable.rvs(alpha=1.5, beta=0, loc=theta, scale=sigma, size=len(returns))
    return levy_process

# Function to plot returns vs Levy process
def plot_returns_vs_levy(returns, levy_process, stock_name, theta, sigma):
    plt.figure(figsize=(10, 5))
    plt.plot(levy_process, label='Generated Levy Process')
    plt.plot(returns.values, label='Actual Returns')
    plt.xlabel('Day')
    plt.ylabel('Returns')
    plt.title(f'{stock_name} Stock Returns vs Generated Levy Process')
    plt.legend()
    plt.show()
    print(f"Estimated Theta (Drift Parameter) for {stock_name}:", theta)
    print(f"Estimated Sigma (Volatility Parameter) for {stock_name}:", sigma)

    # Calculate R-squared
    r_squared = r2_score(returns, levy_process)
    print(f"R-squared for {stock_name}:", r_squared)

# List of stock tickers
stock_tickers = ['AAPL', 'AMZN', 'TSLA']

# Iterate over stocks
for ticker in stock_tickers:
    # Step 1: Retrieve historical stock data
    stock_data = yf.download(ticker, start='2012-01-01', end='2022-12-31')
    returns = stock_data['Adj Close'].pct_change().dropna()

    # Step 2: Estimate parameters using MLE
    theta, sigma = estimate_parameters(returns)

    # Step 3: Generate Levy process using estimated parameters
    levy_process = generate_levy_process(returns, theta, sigma)

    # Step 4: Plot returns vs Levy process
    plot_returns_vs_levy(returns, levy_process, ticker, theta, sigma)

\end{lstlisting}

\end{document}